\documentstyle[12pt]{article}
\newcommand{\const}{\mathop{\rm const}\limits}

\begin{document}

\begin{center}

{\bf  Strengthening of weak convergence for Radon  }\\

\vspace{4mm}

{\bf measures in separable Banach spaces. } \\

\vspace{6mm}

{\sc Ostrovsky E., Sirota L.}\\

\normalsize

\vspace{4mm}

{\it Department of Mathematics and Statistics, Bar-Ilan University,
59200, Ramat Gan, Israel.}\\
e-mail: eugostrovsky@list.ru \\

{\it Department of Mathematics and Statistics, Bar-Ilan University,
59200, Ramat Gan, Israel.}\\
e-mail: sirota3@bezeqint.net \\

\vspace{5mm}

  {\sc Abstract.}

 \vspace{3mm}

 \end{center}

   We prove in this short report that for arbitrary weak converging sequence of sigma-finite Borelian measures in the separable Banach space
there is a compact embedded separable subspace such that this measures not only are concentrated in this subspace but weak converge therein.

\vspace{4mm}

 {\it Key words and phrases:} Separable Banach spaces,  random variables (r.v.), distributions, measures, support, compact
 embedded subspace, Skorokhod's representation, weak convergence, Prokhorov's theorem, distance, sharp quantitative estimate,
 modulus of continuity, Young-Orlicz function, Orlicz spaces.  \par

 \vspace{4mm}

{\it Mathematics Subject Classification (2000):} primary 60G17; \ secondary
 60E07; 60G70.\\

\vspace{3mm}

\section{ Introduction. Notations. Statement of problem. Previous works.}

\vspace{3mm}

 \ Let $  (X = \{x\}, ||\cdot||X)   $   be a separable Banach space relative the norm  function \\
 $ ||\cdot||X; \hspace{5mm}   (\Omega, B, {\bf P}) $ be  a non-trivial  probability space,   $ \xi, \ \{  \xi_n  \}, \ n = 1,2,\ldots  $ be a
 sequence of random variables with values in the space $  X $ having Borel distributions

$$
\mu_n(A) = {\bf P} (\xi_n \in A), \hspace{6mm} \mu(A) = {\bf P} (\xi \in A). \eqno(1.0)
$$

 Recall that the Banach subspace $ (Y = \{  y \}, ||\cdot||Y ) $ of the space $  X $ is named {\it   compact embedded  } into the
space $  X,  $  write

$$
Y \ \stackrel{c.e.}{\subset} \ X,
$$
iff  the space $  Y $ is linear subspace of the space $  X: \ Y \subset X $  and
the closed unit  ball $  B_Y  = \{ y: \ ||y||Y \le 1   \} $ of the space $  Y  $ is pre-compact set in the space $  X,  $ i.e. the closure
$  [B_Y]X $ is compact set  in the space $  X  $ relative the source topology generated by the norm $  ||\cdot||X. $\par

 \ It is known, see \cite{Ostrovsky4}, \cite{Buldygin1}, \cite{Ostrovsky5},  that for one Borelian distribution, say $ \mu = \mu_{\xi},  $
on the space $  X,  $ probabilistic or at last  sigma finite,  there exists a separable compact embedded Banach subspace
$ (Y = \{  y \}, ||\cdot||Y ), $ such that

$$
\mu( X \setminus Y) = 0.
$$

  Note that this proposition is false in the Linear Topological Spaces instead the Banach space  $ X,  $ see \cite{Ostrovsky5}. \par

 Obviously, for the enumerable, or more generally dominated  family  of sigma-finite
 measures $  \mu_n $ there exists a {\it single} Banach separable
compact embedded into the space $  X  $ subspace $ (Y = \{  y \}, ||\cdot||Y ), \ Y \subset X $  such that

$$
\forall n \ \Rightarrow \mu_n( X \setminus Y) = 0. \eqno(1.1)
$$

  Suppose now in addition that the sequence $  \{ \mu_n \} $ converges weakly as $  n \to \infty $ to the measure $  \mu $ in
the classical Prokhorov-Skorokhod sense, i.e.  such that for arbitrary continuous bounded functional $  F: X \to R $

$$
\lim_{n \to \infty} \int_X F(x) \ \mu_n(dx) = \int_X F(x) \ \mu(dx). \eqno(1.2)
$$

 Write $  \mu_n \stackrel{X,w}{\to} \mu. $  \par

 \vspace{3mm}

{\bf Question 1.1.} One can choose either the compact embedded subspace \\
 $ (Y = \{  y \}, ||\cdot||Y ) $  in (1.1) such that the sequence $  \{  \mu_n \}  $ is not only concentrated in the space $  Y,  $
but convergent also in the space  $  Y? $ \par

\vspace{3mm}

{\bf   Our aim this short report is to ground the positive answer on this question.  } \par

\vspace{3mm}

\section{ Main result. }

\vspace{3mm}

{\bf Theorem 2.1.} {\it Let $  (X, ||\cdot||X) $ be  a separable Banach space and $ \mu, \ \mu_n  $ be an enumerable set of Borelian
probability measures (distributions) on $ X $  converging weakly to the measure $ \mu. $ There exists a separable compact embedded into $  X  $
Banach subspace $  (Y, ||\cdot||Y) $ of the space $  X:  \  Y \ \stackrel{c.e.}{\subset} \ X,  $ such that all the measures $ \mu, \ \mu_n  $
are concentrated on the $  Y  $ and moreover weak converge also in the space } $ Y: $

$$
 \mu_n \stackrel{Y,w}{\to} \mu.  \eqno(2.1)
$$

\vspace{4mm}

{\bf Proof.} \par

 {\bf 1.} \ It is sufficient  to consider by virtue of universality only the case when $ X = C[0,1] , $  see
\cite{Ostrovsky4}, \cite{Buldygin1}, \cite{Ostrovsky5}. On the other words, we can and will suppose $ \xi = \xi(t), \ \xi_n =\xi_n(t), \ t \in [0,1]  $
are continuous a.e. numerical values random processes. As  before,

$$
\mu_n \stackrel{C[0,1],w }{\to} \mu.  \eqno(2.2)
$$

\vspace{3mm}

{\bf 2.} Further, we intend to use the famous Skorokhod's representation theorem, \ \cite{Skorokhod1}, see also  \cite{Billingsley1}.
Indeed, there exists a sufficiently rich {\it new probability space} $  (\Omega_1, B_1, {\bf P_1})   $ and identically with $ \xi(t), \ \xi_n(t) $
distributed separable r.p. $  \eta(t), \ \eta_n(t):  $

$$
\xi(t) \stackrel{dis}{=} \eta(t), \hspace{5mm}  \xi_n(t) \stackrel{dis}{=} \eta_n(t), \eqno(2.3)
$$
such that  with the $ {\bf P_1 } $ probability one the sequence $ \{ \eta_n(\cdot) \}   $ converges uniformly to the random process
$  \eta(\cdot). $ \par
 Here the symbol $   \stackrel{dis}{=} $ denotes the coincidence of distribution. Evidently, the r.p. $  \eta, \ \eta_n $
are continuous a.e. \par

 The  corresponding "accompanying" sequence  $  \{ \eta(t), \ \eta_n(t) \} $ is said to be Strengthened Converging Copy of the  sequence
for  the initial one $  \{  \xi(t), \ \xi_n(t) \},   $ not necessary to be unique, write

$$
\{ \eta(t), \ \eta_n(t) \} \stackrel{SCC}{=} \{  \xi(t), \ \xi_n(t) \}.
$$

One can assume

$$
\lim_{n \to \infty} \zeta_n = 0, \hspace{5mm} \zeta_n := \sup_{t \in [0,1]} |\eta_n(t) - \eta(t)|  = 0. \eqno(2.4)
$$

\vspace{3mm}

{\bf 3.} It follows from (2.4) that there exists a deterministic sequence $ \epsilon_n  $ tending to zero and a random variable
$  \tau  $  defined on the new probability space such that

$$
\zeta_n \le \tau \cdot \epsilon_n, \eqno(2.5)
$$
see e.g. \cite{Kantorovich1}, chapter 2, section 3. \par
 Moreover, since the sequence  of continuous functions $ \eta(t), \ \eta_n(t) $ converges uniformly  (a.e.),
therefore it is compact set in the space $  C[0,1]. $  It follows from the Arzela-Ascoli theorem that they are  equicontinuous, i.e. there
exists the (random) non-negative  continuous increasing function $ \delta \to h(\omega, \delta), \ \delta \in [0,1]: $

$$
\lim_{\delta \to 0+} h(\omega, \delta) = 0, \eqno(2.6)
$$
such that

$$
\Delta(\eta_n - \eta, \delta) \le h(\omega, \delta) \to 0, \ \delta \to 0+. \eqno(2.7)
$$

 \ In what follows $ \Delta(f, \delta)  $ will be denote the ordinary module (modulus) of continuity of an uniform continuous
function $ f \in C[0,1].  $  \par

\vspace{3mm}

 {\bf 4.} It follows from the main result of the preprint  \cite{Ostrovsky4} that there exists a random variable $   \theta $
defined on the new probability space and {\it deterministic } non-negative continuous function $ \delta \to g(\delta),  $
which takes zero value at the origin:  $  g(0) = g(0+) = 0  $ such that

$$
h(\omega, \delta) \le \theta \cdot g(\delta). \eqno(2.8)
$$
 Then

$$
\Delta(\eta_n - \eta, \delta) \le \theta \cdot g(\delta). \eqno(2.9)
$$

\vspace{3mm}

{\bf 5.} Let us introduce the following modification  of the classical H\"older's spaces $ H^o(\sqrt{g}). $  By definition,
 $ H^o(\sqrt{g}) $  consists on all the (continuous) functions $ f: [0,1] \to R $ for which

$$
\lim_{\delta \to 0+} \frac{\Delta(f,\delta)}{\sqrt{ g(\delta)  }} = 0, \eqno(2.10)
$$
with (finite) norm

$$
||f||H^o(\sqrt{g}) \stackrel{def}{=} \max_{t \in [0,1]} |f(t)| + \sup_{\delta \in (0,1)} \left[ \frac{\Delta(f,\delta)}{\sqrt{ g(\delta)  }} \right]. \eqno(2.11)
$$
 These Banach spaces are separable and compact embedded into the space $ C[0,1], $  see, e.g. the monograph
 \cite{Guseinov1}, chapter 1, where these spaces are used in particular in the theory of non-linear singular integral equation;
 some another applications, for example, in the theory of CLT in Banach spaces, may be found in \cite{Ostrovsky9}.\par

\vspace{3mm}

{\bf 6.}  We can  suppose without loss of generality that the r.f. $  \eta(\cdot), \ \eta_n(\cdot) $ belong to the introduced above  space
$  H^o(\sqrt{g}). $ \par

 Indeed, as long as the sequence $  \eta_n(\cdot) $ converges uniformly  to $ \eta(\cdot), $  it is also compact.  It remains to repeat the
considerations of fourth item (2.9):

$$
\Delta(\eta_n, \delta) \le \tilde{\theta} \cdot \tilde{g}(\delta),
$$

$$
\Delta(\eta, \delta) \le \tilde{\theta} \cdot \tilde{g}(\delta),
$$
and choose $ \hat{g}(\delta) := \max(\tilde{g}(\delta), \ g(\delta)).  $ \par

\vspace{4mm}

{\bf 7.} It follows immediately from the equality (2.9) that the sequence of the r.p. $  \{  \eta_n(\cdot) \} $ converge as $  n \to \infty $
to the r.p. $  \eta(\cdot) $ also in the norm $  H^o(\sqrt{g}) $ with $ {\bf P_1}  $ probability one:

$$
|| \eta_n(\cdot) - \eta(\cdot)|| H^o(\sqrt{d}) \to 0. \eqno(2.12)
$$

 The last equality implies for the source sequence of r.p. $  \xi(\cdot), \ \xi_n(\cdot) $ weak its distribution convergence in the space
$ H^o(\sqrt{d}). $  \par

 This completes the proof of theorem 2.1.\par

\vspace{4mm}

{\bf Example 2.1.} Suppose that the sequence of centered continuous random fields $  \xi_n(s), \ s \in S, $ somehow dependent,
where $  (S = \{ s \},r) $ is compact relative certain distance $ r(\cdot, \cdot) $ metric space, converges weakly  to the continuous
{\it Gaussian} random field $  \xi = \xi(s) $  in the ordinary space of all continuous functions $  C(S); $  on the other words,
CLT in $ C(S), $ \cite{Ledoux1}, chapter 9, \  \cite{Ostrovsky1}, chapter 4; uniform CLT \cite{Dudley1}. \par

  We deduce based on the theorem 2.1 that there exists some modified H\"older space $ H^o(g)  $  over $  (S,r) $
such that the sequence $  \xi_n(\cdot) $  convergent weakly in the space $ H^o(g), $ i.e. is subgaussian,  as well. \par

\vspace{3mm}

\section{ Orlicz norm estimates for the tail of random coefficient.} \par

\vspace{3mm}

{\bf  The case of  the space of continuous functions.}\par

\vspace{3mm}

 \ In this subsection $  X = C(S),  $ where the set $ S = \{ s \}  $ is compact set relative certain distance $  \rho = \rho(s_1, s_2),
 \ s_1, s_2 \in S. $ \par

 \ It is interest by our opinion for the practical using, for instance, in the Monte - Carlo method, to estimate the tails of distribution
for the r.v. $  \theta $ in the estimation 2.9.  For this purpose assume that the r.v.

$$
\nu \stackrel{def}{=} \sup_n ||\xi_n||C(S) = \sup_n \sup_{s \in S} |\xi_n(s)| \eqno(3.1)
$$
belongs to certain Orlicz space $  L(\Phi)  $ with Luxemburg norm  $ || \cdot ||L(\Phi) = || \cdot ||\Phi $ constructed over source probability
space; here $  \Phi(\cdot) $ is some Young-Orlicz function.  \par

 The detail investigation of Orlicz's spaces may be found in the classical books \cite{Krasnoselsky1},  \cite{Rao1}, \cite{Rao2}.

 We can and will suppose without loss of generality

$$
{\bf E} \Phi(\nu) = 1. \eqno(3.2)
$$

\vspace{3mm}

{\bf Proposition 3.1.} \par

\vspace{3mm}

 \ It follows in particular from one of results of the recent  preprint \cite{Ostrovsky10}
that if the function $ \Phi $ satisfies the so-called $   \Delta_2 $ condition:  $ \Phi \in \Delta_2,  $ then
the {\it scaling} function $  g = g(\delta) $ in (2.9) may be picked such that also  $ \theta  \in L(\Phi). $ \par

 The converse predicate is trivially true. \par

 In the opposite case the situation is more complicated. Recall that the other Orlicz function $   \Psi(\cdot)  $ is called weaker than the
function $  \Phi, $  notation $   \Psi << \Phi, $ if for all positive constant $  v; v = \const > 0 $

$$
\lim_{u \to \infty} \frac{\Psi(u v)}{\Phi(u)} = 0. \eqno(3.3)
$$

 It is alleged that  if the function $ \ \Psi, \Psi << \Phi  $ is given, then
 the {\it scaling} function $  g = g(\delta) = g_{\Phi, \Psi}(\delta) $ may be picked such that also  $ \theta  \in L(\Psi). $ \par

\vspace{4mm}

{\bf Example 3.1.}   For instance, if

\vspace{4mm}

$$
 \sup_n {\bf E} | \ ||\xi_n||C(S) \ |^p < \infty, \ \exists p = \const \ge 1,
$$
then  the  scaling function $  g = g(\delta) $ may be picked such that also  $ {\bf E} | \theta|^p  = 1. $\par

\vspace{4mm}

{\bf Example 3.2.}

\vspace{4mm}

  Let us consider now  the so-called Gaussian centered case,  i.e. when the {\it  common } distribution of the infinite-dimensional vector
$  \vec{\xi} = \{  \xi, \xi_1, \xi_2, \ \xi_3,  \ldots   \} $   has a  mean zero Gaussian distribution.  We have to take

$$
\Phi(u) = \Phi_G(u) := e^{u^2/2} - 1.
$$

 \ It follows from one of the main results of an article X.Fernique \cite{Fernique1}, see also \cite{Buldygin1}, \cite{Ostrovsky10} that
for any choice of the function $  g(\delta) $ satisfying the relation (2.9) the r.v. $ \theta $ belongs also to the Orlicz's space
$   L(\Phi_G), $ i.e. is subgaussian. \par

 \vspace{3mm}

  The last example show us that the second assertion of proposition 3.1  is in general case improvable.\par

\vspace{4mm}

 \section{ Bernstein's moment convergence for compact embedded subspace. }

\vspace{3mm}

{\bf The general case of  the arbitrary separable Banach space. } \par

\vspace{3mm}

 \ The Bernstein's moment convergence for  weakly convergent sequence of measures $ \{  \mu_n \}   $ imply by definition
the following integral convergence

$$
\int_X V(x) \ \mu_n(dx)  \to \int_X V(x) \ \mu(dx)
$$
for certain continuous {\it unbounded} functional $ V: X \to R.  $ This problem goes back to  S.N.Bernstein   \cite{Bernstein1}; see also
\cite{Hall1}, \cite{Peligrad1}. \par

 As a rule, in the aforementioned articles the functional $  V(x) $ has a form $  V(x) = ||x||^p, \ x \in X, \ p = \const \ge 2.  $\par

\vspace{3mm}

 \ Let again $ (X, \ ||\cdot||X ) $ be separable Banach space and let $ \vec{\mu} = ( \mu, \  \{  \mu_n \} ) $  be the family of Borelian probability
measures defined on all the Borelian subsets $  X.  $ Let also $  V: X \to R $  be continuous functional acting  from $  X  $ to the real  axis
$  R: \ V = V(x), \ x \in X.  $ \par

\vspace{3mm}

 \ {\bf Definition 4.1.}  The functional $  V(\cdot) $ is named
 uniform integrable relative the family $ \vec{\mu}, $ if

 $$
 \lim_{N \to \infty} \sup_n \int_{x: |V(x)| > N } |V(x)| \ \mu_n(dx) = 0. \eqno(4.1)
 $$

\vspace{3mm}

{\bf Lemma 4.1.} {\it Suppose that $  \mu_n \stackrel{w,X}{\to} \mu  $ and that the functional $  V(\cdot) $ is  uniform integrable
relative the family $ \vec{\mu}, \ \mu. $  We propose}

$$
\lim_{n \to \infty} \int_X V(x) \ \mu_n(dx) = \int_X V(x) \ \mu(dx). \eqno(4.2)
$$

\vspace{3mm}

{\bf Proof } is elementary.  Let the positive number $  \epsilon > 0 $ be a given.  We  introduce the truncated functional,
also continuous,

$$
V_N(x) = V(x), \ |V(x)| \le N; \hspace{6mm} V_N(x) = - N, \ V(x) < - N;
$$

$$
 V_N(x) = + N, \ V(x) >  N; \hspace{6mm} N = 1,2,\ldots,
$$
and denote

$$
\kappa := \left| \ \int_X V(x) \ \mu_n(dx) - \int_X V(x) \ \mu(dx) \ \right|.
$$

 We get using the triangle inequality  $  \kappa \le \kappa_1 + \kappa_2 + \kappa_3, $ where

$$
\kappa_1 =  \left| \int_X V(x) \ d \mu_n -   \int_X V_N(x) \ d \mu_n  \right|,
$$

$$
\kappa_2 =  \left| \int_X V_N(x) \ d \mu_n -   \int_X V_N(x) \ d \mu   \right|,
$$

$$
\kappa_3 =  \left| \int_X V(x) \ d \mu -   \int_X V_N(x) \ d \mu  \right|.
$$

 Since the functional $  V = V(x)  $ is  uniform integrable relative the family $ \vec{\mu}, $ there exists the value $  N_0 = N_0(\epsilon)  $
such that  for all the values $  N > N_0(\epsilon)   $

$$
\kappa_1 < \epsilon/3, \ \hspace{6mm} \kappa_3 \le \epsilon/3.
$$
 Further, as long as $ \mu_n \stackrel{X,w}{\to} \mu, $ there is a value $ n_0 = n_0(\epsilon, N_0(\epsilon))  $ so that
for all the values  $  n > n_0 $

$$
\kappa_2 < \epsilon/3.
$$
 Totally, $  n > n_0 \ \Rightarrow  \kappa < \epsilon, $ Q.E.D. \par

 \vspace{3mm}

 We deduce  applying this assertion and the last section the following proposition. \par

 \vspace{3mm}

{\bf Theorem 4.1.}  {\it Let all the conditions of Lemma 4.1 be satisfied. There exists common support for the
all the measures $ \mu, \ \mu_n  $ compact embedded Banach subspace
$  (Y, \ ||\cdot||Y) $  such that the functional $  V(\cdot)  $ is uniform integrable also for these measures
inside the new space $  Y $ and hence}

$$
\lim_{n \to \infty} \int_Y V(y) \mu_n(dy) =  \int_Y V(y) \mu(dy). \eqno(4.3)
$$

\vspace{4mm}

\section{ Concluding remarks. }

\vspace{3mm}

\ {\bf A. Generalization  on the non-normed measures.} \par

\vspace{3mm}

 \ All the references about the random variables $ \ \xi, \ \xi_n  $  generating  the correspondent distributions
 $ \  \mu, \ \mu_n  $ may be eliminated; it is sufficient to consider only
the sequence of Borelian (Radon) sigma-finite measures $  \{ \mu_n \} $ converging weakly in at  the same separable Banach space $ X $
 to the  measure $  \mu, $ which also also Borelian.\par
 Wherein theorem 2.1 remains true, as long as each Borelian sigma-finite measure in equivalent in the Radon-Nikodym sense to the
probability distribution.  \par

 \vspace{3mm}

\vspace{3mm}

\ {\bf B. \ Generalization on the arbitrary family of  measures.} \par

\vspace{3mm}

 \ At the same result (theorem 2.1) remains true   for an arbitrary  dominated and convergent {\it net} $  \mu_{\alpha}, \ \alpha \in A,
 \ \alpha \to \alpha_0 \in A $
 of sigma-finite  Borelian measures in the separable Banach space  $ X. $ \par

\vspace{3mm}

 \ {\bf C. Open question.} \par

\vspace{3mm}

 V.V.Buldygin in \cite{Buldygin1} proved that in the probabilistic case $ \ \mu(X) = \mu_n(X) = 1 \ $ the single subspace
$  (Y, ||\cdot||) $ may be constructed to be reflexive and with continuous differentiable on the unit sphere
in the Freshet sense norm. \par

  We do not know either or not possible to  choose this "good" subspace $  (Y, ||\cdot||) $  such that all the r.v.
are concentrate in this space and in addition weakly converges therein.\par

\end{document}